\documentclass[12pt,leqno]{article}

\usepackage{amssymb, graphicx, amscd, amsmath, amssymb, amsthm, array, fancyvrb, relsize, paralist,url}

\numberwithin{equation}{section}

\def\Fo{\mathbb{F}}
\def\F{\Fo_p}
\def\fx4{\frac{x}{4}}

\def\3F2{{}_3\hspace{-1pt}F_2}
\def\2F1{{}_2\hspace{-1pt}F_1}
\def\h2F1{{}_2\hspace{-1pt}\widehat{F}_1}
\def\Fq2{\F_{q^2}}

\def\ep{\epsilon}
\def\c4{\chi_4}
\def\oc4{\overline{\chi_4}}

\def\l({\left(}
\def\r){\right)}
\def\bar{\begin{array}{r|}}
\def\ear{\end{array}}

\def\CycFull{\mathrm{CycFull}}
\def\CycShort{\mathrm{CycShort6}}
\def\CycShorter{\mathrm{CycShort8}}

\def\ind{\mathrm{ind}}

\title{Nonexistence of twenty-fourth power residue addition sets}

\author{\\ \\ Ron Evans\\
Department of Mathematics\\
University of California at San Diego\\
La Jolla, CA  92093-0112 \\
revans@ucsd.edu
\\ \\
and
\\ \\
Mark Van Veen\\
Varasco LLC \\ 
2138 Edinburg Avenue\\ 
Cardiff by the Sea, CA 92007  \\
mark@varasco.com
}

\medskip

\date{November 2016}

\begin{document}

\maketitle


\newpage

\begin{abstract}
Let $n>1$ be an integer, and let $\F$ denote a field of $p$ elements 
for a prime $p \equiv 1 \pmod n$.
By 2015, the question of existence or nonexistence
of $n$-th power residue difference sets in $\F$ 
had been settled for all $n <24$.
We settle the case $n=24$ by proving
the nonexistence of $24$-th power residue difference sets in $\F$.
We also prove
the nonexistence of \textit{qualified} $24$-th power residue difference sets
in $\F$.
The proofs make use of a Mathematica program which computes 
formulas for the cyclotomic numbers of order 24 in terms of parameters
occurring in quadratic partitions of $p$.
\end{abstract}

\maketitle

\section{Introduction}
For an integer $n >1$, let $p$ be a prime of the form $p=nf+1$.
Let $H_n$ denote the set of nonzero $n$-th power residues in $\F$,
where $\F$ is the field of $p$ elements.  For $\ep \in \{0,1\}$,
define $H_{n, \ep} = H_n \cup \{1-\ep\}$. Note that the set
$H_{n, \ep}$ has $f+\ep$ elements.

Fix $m \in \F^*$.  Lam \cite{Lam} called $H_{n, \ep}$ an $n$-th power
residue addition set if the list of differences $s - mt \in \F^*$
with $s,t \in H_{n, \ep}$ hits each element of $\F^*$ the same number of
times. If $m$ is an $n$-th power residue,
such an addition set is called an $n$-th power
residue difference set.  If $m$ is not an $n$-th power residue,
then as in \cite[p. 94]{Byard},
such an addition set is called a qualified $n$-th power
residue difference set with qualifier $m$.

Let $g$ denote a primitive root modulo $p$.   For integers $s,t$ modulo $n$,
the cyclotomic number $C_n(s,t)$ of order $n$ is defined to be the number
of integers $N \in \F$ for which both $N/g^s$ and $(N+1)/g^t$
are nonzero $n$-th power residues in $\F$.
If $H_{n, \ep}$ is a difference set, then necessarily \cite[p. 677]{Ev24}
$n$ is even, $f$ is odd, and
\begin{equation}\label{1.1}
n^2C_n(s,0) = p-1+2n\ep -n, \quad 1 \le s < n/2.
\end{equation}
If $H_{n, \ep}$ is a qualified difference set, then necessarily 
\cite[Theorems 2.1 and 2.2]{Byard} $n$ is even, $f$ is even, and
\begin{equation}\label{1.2}
n^2C_n(s,n/2) = p-1+2n\ep, \quad 1 \le s < n/2.
\end{equation}

For $n < 24$, it is known that 
$H_{n, \ep}$ can be a difference set only in the three exceptional cases
$H_{2, \ep}$, $H_{4, \ep}$, $H_{8, \ep}$ 
listed in \cite[(1.1)--(1.3)]{Byard}.
This follows from the work of a number of different authors during
the period 1933--2015. For references, consult Xia \cite{Xia}, who
has extended the results to fields of $q$ elements, where $q$ is a prime
power.
In Section 4, we prove that $H_{24, \ep}$ cannot be a difference set,
by showing that \eqref{1.1} cannot hold for $n=24$.

For $n < 22$, it is known that 
$H_{n, \ep}$ can be a 
qualified difference set only in the three exceptional cases
$H_{2, \ep}$, $H_{4, \ep}$, $H_{6, \ep}$ 
listed in \cite[(1.4)--(1.6)]{Byard}.
In Section 3, we prove that $H_{24, \ep}$ cannot be a 
qualified difference set,
by showing that \eqref{1.2} cannot hold for $n=24$.

Our proofs depend on formulas for the cyclotomic numbers $C_{24}(s,t)$.
Printed tables of these formulas were archived in 1979 \cite{EvUMT}, but it is
much more useful to have digital access.  Thus we wrote a Mathematica
program \cite{EvMat} to compute the formulas for $C_{24}(s,t)$.   This program
is described in the next section.

We remark that besides their use for analyzing power residue difference
sets, cyclotomic numbers have applications to such topics as
counting points on elliptic curves \cite{XiaYang}, Gauss periods and 
complexity of normal bases for finite fields \cite{CGPT},\cite{GaoT},
cyclic codes \cite{Ding}, cryptographic functions \cite{CDR},
residuacity \cite[Chapter 7]{BEW}, linear complexity of sequences
\cite{BrandW}, and almost difference sets \cite{DPW}.

\section{Cyclotomic numbers of order 24}

Let $\beta = \exp(2 \pi i/24)$.
For $0 \le u,v \le 23$, define
the Jacobi sum $J(u,v,\beta)$ by
\begin{equation}\label{2.1}
J(u,v)=J(u,v,\beta) = \sum_{x=2}^{p-1}\beta^{\ind(x)u + \ind(1-x)v},
\end{equation}
where $\ind(x)$ denotes the index of $x$ with respect to the primitive
root $g$.  The cyclotomic numbers $C_{24}(s,t)$ can be computed in
terms of Jacobi sums via the formula \cite[eq. 2.5.1]{BEW}
\begin{equation}\label{2.2}
576C_{24}(s,t)=\sum_{u=0}^{23} \sum_{v=0}^{23} 
(-1)^{uf} \beta^{-su - tv} J(u,v,\beta).
\end{equation}

There are five Jacobi sums in \eqref{2.2}
that have been  expressed in \cite[p. 678]{Ev24} in terms of 
sixteen integer parameters called
\begin{equation}\label{2.3}
X, Y, A, B, C, D, U, V, D_j, \quad 0 \le j \le 7,
\end{equation}
viz.
\begin{equation}\label{2.4}
J(6,12)=-X+2Yi, \quad (p=X^2+4Y^2, \quad X\equiv 1 \pmod 4),
\end{equation}
\begin{equation}\label{2.5}
J(4,12)=-A+Bi\sqrt{3}, \quad (p=A^2+3B^2, \quad A\equiv 1 \pmod 6),
\end{equation}
\begin{equation}\label{2.6}
J(3,12)=-C+Di\sqrt{2}, \quad (p=C^2+2D^2, \quad C\equiv 1 \pmod 4),
\end{equation}
\begin{equation}\label{2.7}
J(1,12)=U+2Vi\sqrt{6}, \quad (p=U^2+24V^2, \quad U\equiv -C \pmod 3),
\end{equation}
\begin{equation}\label{2.8}
J(1,2)=\sum_{j=0}^7 D_j \beta^j.
\end{equation}
The remaining Jacobi sums in \eqref{2.2} are expressible
\cite[Chapter 3]{BEW} in terms of the parameters \eqref{2.3}
together with the four parameters
\begin{equation}\label{2.9}
F1, V1, Z, T,
\end{equation}
where $F1 \in \{0,1\}$ with $F1 \equiv f \pmod 2$,
$V1 \in \{0,1\}$ with $V1 \equiv V \pmod 2$,
$Z \equiv \ind(2) \pmod {12}$, and
$T \equiv \ind(3) \pmod 8$.
As described in \cite[p. 678]{Ev24}, $g$ may be chosen so that
$Z \in \{0,2,4,6\}$ and $T \in \{0,2,4\}$.
Thus there are 48 distinct 4-tuples $\{F1, V1, Z, T\}$, and for each
such 4-tuple, one can create a table of the 576 numbers
\begin{equation}\label{2.10}
576C_{24}(s,t), \quad 0 \le s, t \le 23.
\end{equation}
Each number in \eqref{2.10} turns out to be an integer  linear combination of
\begin{equation}\label{2.11}
p, 1, X, Y, A, B, C, D, U, V, D_0, D_1, D_2, D_3, D_4, D_5, D_6, D_7.
\end{equation}
For example, when $\{F1, V1, Z, T\}=\{1,1, 4,0\}$, we have
\begin{equation}\label{2.12}
\begin{split}
& 576C_{24}(6,0) =  p-23+4X+0Y-14A+24B-8C+0D -8U \\
&+0V+32D_0+0D_1+0D_2+0D_3+16D_4+0D_5+0D_6+0D_7.
\end{split}
\end{equation}

Our Mathematica module $\CycFull[s,t,F1, V1, Z, T]$ in \cite{EvMat}
outputs a list beginning with $s,t$ followed in order by the coefficients
in the linear combination of the elements in \eqref{2.11}.   For example,
in view of \eqref{2.12}, the input $\CycFull[6,0,1,1, 4, 0]$ outputs the list
\[
\{6,0,1,-23,4,0,-14,24,-8,0,-8,0,32,0,0,0,16,0,0,0\}.
\]
For each of the 48 tuples $\{F1, V1, Z, T\}$, CycFull can be used
to construct a table of the 576 cyclotomic numbers of order 24.

We now describe two additional Mathematica modules, 
CycShort6 and CycShort8,
which will be used in the sequel.
CycShort6 is like CycFull except that the output list is shortened
by omitting the first four entries and also omitting the coefficients
of the six parameters $D, V, D_1, D_3, D_5, D_7$.
CycShort8 is like CycShort6 except that the output list is further
shortened by omitting the coefficients of the two parameters
$D_2, D_6$.   For example, $\CycShort[6,0,1,1,4,0]$ outputs
the list
\[
\{4,0,-14,24,-8,-8,32,0,16,0\}
\]
of coefficients of
\begin{equation}\label{2.13}
X,Y,A,B,C,U,D_0,D_2,D_4,D_6,
\end{equation}
while $\CycShorter[6,0,1,1,4,0]$ outputs
the list
\[
\{4,0,-14,24,-8,-8,32,16\}
\]
of coefficients of
\begin{equation}\label{2.14}
X,Y,A,B,C,U,D_0,D_4.
\end{equation}

\section{Nonexistence of qualified difference sets}

In this section we focus on the 24 tuples $\{0,V1,Z,T\}$, which correspond
to Tables 1--24 in \cite{EvMat}.  Running  $\CycFull[s,12,0,V1,Z,T]$,
we see that
\begin{equation}\label{3.1}
576C_{24}(s,12) = p+1 + \gamma_0(s), \quad 1 \le s \le 11,
\end{equation}
where
\[
\gamma_0(s)=d_1X+d_2Y +d_3A+d_4B+d_5C+d_6U+d_7D_0+d_8D_2+d_9D_4+d_{10}D_6
\]
for integer coefficients $d_j$ depending on $s, V1, Z, T$.
Thus the coefficients of $\gamma_0(s)$ are given by $\CycShort[s,12,0,V1,Z,T]$.
Moreover, in the eight cases where $T=4$, the coefficients $d_8$
and $d_{10}$ are always 0, so in those eight cases,
the coefficients of $\gamma_0(s)$ are given by $\CycShorter[s,12,0,V1,Z,T]$.

Assume for the purpose of contradiction that \eqref{1.2} holds.
Then by \eqref{3.1}, we have the following system of eleven linear equations:
\begin{equation}\label{3.2}
\gamma_0(s)=-2+48\ep, \quad 1 \le s \le 11.
\end{equation}
We will  prove that the system \eqref{3.2}
has no viable solution, thus obtaining the desired result
that \eqref{1.2} cannot hold.

The system \eqref{3.2} can be represented as a matrix equation of the form
\begin{equation}\label{3.3}
M\boldsymbol{y}=\boldsymbol{h},
\end{equation}
where $M$ is an 11 by 10 matrix of integer coefficients,
$\boldsymbol{y}$ is the column vector whose ten entries are the variables
in \eqref{2.13}, and $\boldsymbol{h}$ is the column vector whose eleven
entries all equal $-2+48\ep$.  In the eight cases where $T=4$,
we can instead take $M$ as an 11 by 8 matrix of integer coefficients
with $\boldsymbol{y}$ the column vector whose eight entries are the variables
in \eqref{2.14}. Our  matrix equations are displayed in 
\cite{EvMat} in Tables 1--24.
In several cases, our matrix equations don't need to 
make use of all eleven rows of the
matrix $M$. For example,
for Table 21, we obtain a contradiction with a coefficient 
matrix consisting  only
of rows 4 through 11, i.e., the top three rows of $M$ have been omitted.

For each of the 24 tables, a particular solution to the matrix equation 
is found using the Mathematica function LinearSolve.
To find the general solution, we add the particular solution to the 
null space of the matrix, computed in the tables with the function NullSpace.
We proceed to examine the 
general solutions table by table, and show
that none of them are viable.

\noindent
\section*{Table 1} 

As shown in \cite{EvMat},
in the general solution to the matrix equation, we have
\[
X=A=1-24\ep, \ B=2b, \ D_0= -1+24\ep-b
\]
for some $b$, and $b$ must be an integer since $D_0$ is.
By the formulas for $p$ in \eqref{2.4}--\eqref{2.5},
we have $4Y^2 = 12b^2$.   This contradicts the fact that
3 is not a square, so the matrix equation has no viable solution.

\noindent
\section*{Tables 2,3,7,9,10,11,12,13,14,15,23,24}

In the general solution, $B=0$, which contradicts \eqref{2.5}.

\noindent
\section*{Tables 4,5,6,16,17,18,19,20,21}

In the general solution, $X$ is not an integer, 
which is a contradiction.

\noindent
\section*{Table 8}

In the general solution,
\[
X=1-24\ep-16a, \ Y=a, \ A=-1+24\ep+8a, \ B=-4a
\]
for some $a$, and $a$ must be a nonzero integer since $Y$ is.
By the formulas for $p$ in \eqref{2.4}--\eqref{2.5},
we have 
\[
p=(1-24\ep-16a)^2+4a^2=(-1+24\ep+8a)^2+48a^2.
\]
Solving for $a$ yields $a=(4-96\ep)/37$, which contradicts
the fact that $a$ is an integer.

\noindent
\section*{Table 22}

In the general solution, for some $a$, we have
\[
X=(-1+24\ep)/23 + 192a, \ Y=92a, \ A=9(-1+24\ep)/23+256a, \ B=138a.
\]
By the formulas for $p$ in \eqref{2.4}--\eqref{2.5}, we have
\[
((-1+24\ep)/23 + 192a)^2 +4(92a)^2=(9(-1+24\ep)/23+256a)^2+3(138a)^2.
\]
Solving for $a$ yields $a=(2-48\ep)/897$ or $a=(10-240\ep)/7659$.
In particular, $46a$ is not an integer.
However, in the general solution, we also have $46a = D_7-D_6$,
which is not possible since $D_7$ and $D_6$ are integers.

\vspace{1 mm}

This completes the proof that 
$H_{24, \ep}$ cannot be a
qualified difference set.

\section{Nonexistence of difference sets}

In this section we focus on the 24 tuples $\{1,V1,Z,T\}$, which correspond
to Tables 25--48 in \cite{EvMat}.  Running  $\CycFull[s,0,1,V1,Z,T]$,
we see that
\begin{equation}\label{4.1}
576C_{24}(s,0) = p -23 + \gamma_1(s), \quad 1 \le s \le 11,
\end{equation}
where
\[
\gamma_1(s)=c_1X+c_2Y +c_3A+c_4B+c_5C+c_6U+c_7D_0+c_8D_2+c_9D_4+c_{10}D_6
\]
for integer coefficients $c_j$ depending on $s, V1, Z, T$.
Thus the coefficients of $\gamma_1(s)$ are given by $\CycShort[s,0,1,V1,Z,T]$.
Moreover, in the eight cases where $T=0$, the coefficients $c_8$
and $c_{10}$ are always 0, so in those eight cases,
the coefficients of $\gamma_1(s)$ are given by $\CycShorter[s,0,1,V1,Z,T]$.

Assume for the purpose of contradiction that \eqref{1.1} holds.
Then by \eqref{4.1}, we have the following system of eleven linear equations:
\begin{equation}\label{4.2}
\gamma_1(s)=-2+48\ep, \quad 1 \le s \le 11.
\end{equation}
In the same manner used in the previous section, 
we will  prove that the system \eqref{4.2}
has no viable solution, thus obtaining the desired result
that \eqref{1.1} cannot hold.

\section*{Tables 25,26,27,28,30,31,33,34,36,38,39,47,48}

In the general solution, we have $B=0$ or $Y=0$,
both of which are impossible.

\section*{Tables 29,32,35,41,42,43,44,45,46}

In the general solution, $X$ is not an integer, 
which is a contradiction.

\section*{Table 37}

In the general solution,
\[
X=5-120\ep, \ A=13-312\ep, \ C=-23+552\ep, \ U=-1+24\ep
\]
and for some $a$,
\[
B=-2a, \ D_0 = 8 - 192\ep -a, \ D_4=2a.
\]
Here $a$ must be a nonzero integer since $D_0$ is an integer and $B$ is nonzero.
We do not see how to obtain a contradiction via the systematic method
used for the other tables.   However, a contradiction has been obtained
using Gauss sums of order 24;  see \cite[pp. 680--683]{Ev24}.
(\textit{Note:} On lines 22, 26, 30, 32 of \cite[p. 680]{Ev24},
subtract 7 from each reference number.  Also, on line 2 of the Introduction
in \cite{Ev24}, replace the second $p$ by $e$.)

\section*{Table 40}

In the general solution, for some $a$,
\begin{equation*}
\begin{split}
&X=(-5+120\ep)/19 - 24a, \ Y=-19a, \\
&A=(-17+408\ep)/19-112a, \ B=(7-168\ep)/19+26a.
\end{split}
\end{equation*}
Here $19a$ is an integer, since $Y$ is.
By the formulas for $p$ in \eqref{2.4}--\eqref{2.5}, we have
\begin{equation*}
\begin{split}
&((-5+120\ep)/19 - 24a)^2 +1444a^2 =\\
&((-17+408\ep)/19-112a)^2+3((7-168\ep)/19+26a)^2.
\end{split}
\end{equation*}
This equation has irrational solutions $a$, which is a contradiction.

\vspace{1 mm}

This completes the proof that
$H_{24, \ep}$ cannot be a
difference set.

\end{document}